\documentclass{amsart}

\usepackage[all]{xy} 
\input xy 

\xyoption{all}
\xyoption{2cell} 
\xyoption{v2}

\DeclareMathOperator{\Hom}{Hom} 

\theoremstyle{plain}
\newtheorem{theorem}{Theorem}[section]

\newtheorem{proposition}[theorem]{Proposition}

\theoremstyle{definition}
\newtheorem{definition}[theorem]{Definition}

\theoremstyle{remark}

\newtheorem{note}{Note}[section]

\numberwithin{equation}{section}
\numberwithin{figure}{section}

\newcommand{\cU}{{\mathcal U}}

\newcommand{\cV}{{\mathcal V}}

\newcommand{\RR}{{\mathbb R}}
\newcommand{\ZZ}{{\mathbb Z}}
\newcommand{\QQ}{{\mathbb Q}}

\renewcommand{\a}{\alpha}
\renewcommand{\b}{\beta}
\renewcommand{\c}{\gamma}

\begin{document}

\title[Bundle gerbes and infinite-dimensionality]{A note on bundle gerbes and infinite-dimensionality}
  \author[M. Murray]{Michael Murray}
  \address[Michael Murray]
  {School of Mathematical Sciences\\
  University of Adelaide\\
  Adelaide, SA 5005 \\
  Australia}
  \email{michael.murray@adelaide.edu.au}

\author[D. Stevenson]{Danny Stevenson}
\address[Danny Stevenson]
{Department of Mathematics\\
University of Glasgow  \\
15 University Gardens\\
Glasgow G12 8QW\\
United Kingdom}
\email{d.stevenson@maths.gla.ac.uk} 
   
  \thanks{The first author acknowledges the support of the Australian
Research Council. Both authors thank Diarmuid Crowley for explaining  Hurewicz's theorem to them.}


\begin{abstract}
Let $(P, Y)$ be a bundle gerbe over a fibre bundle $Y \to M$.  We show that if $M$ is simply-connected and
the fibres of $Y \to M$  are connected and finite-dimensional then the Dixmier-Douady class of $(P, Y)$ is torsion.
This corrects and extends an earlier result of the first author. 
\end{abstract}

\maketitle

\section{Introduction}

This paper is dedicated to Alan Carey in commemoration of his 60th birthday.  The idea of bundle gerbes \cite{Mur96Bundle-gerbes} had its original  motivation in  attempts  by the 
first author and Alan to geometrize degree three cohomology classes. This in turn arose from a shared interest  in anomalies in quantum field theory resulting from non-trivial cohomology 
classes in the space of connections modulo gauge. Even in the earliest of our joint papers on anomalies  \cite{CarMur}, 
which demonstrates that the Wess-Zumino-Witten term can be understood as holonomy for a line bundle on the loop group, there is a bundle gerbe lurking, at that time unnoticed, in the background.  It was not till sometime later that we realised  that a better interpretation of the Wess-Zumino-Witten term for a map of a surface into a compact Lie group is as the surface holonomy
of the pull-back of the basic bundle gerbe over that group \cite{CarMicMur}. 

In the present work we are concerned with the relationship between bundle gerbes and infinite-dimensionality.
It is well-known \cite{AtiSeg, BouCarMat} that there is a distinct difference in twisted $K$-theory
over a manifold $M$  between the case where the twist $\a \in H^3(M, \ZZ)$ is torsion and the case where it is  of infinite order. The latter
seems to necessitate infinite-dimensional constructions in a way that the former does not. A similar
situation holds in the case of geometric realisations of the twist $\a$ as gerbes and bundle gerbes. In particular
in \cite{Mur96Bundle-gerbes} it was claimed by the first author that the following was true.

\begin{theorem}
\label{thm} 
Let $Y \to M$ be a fibre bundle with finite-dimensional $1$-connected fibres. Let $M $ also be $1$-connected. 
Then any bundle gerbe $(P, Y)$ over $M$ has exact three-curvature and hence torsion Dixmier-Douady class.
\end{theorem}

Unfortunately the proof given in \cite{Mur96Bundle-gerbes} is incorrect. We will explain why this is the case and give a correct proof below.  Moreover we will extend this result to the case that the fibre is just connected. 
In addition we will give examples of bundle gerbes with non-torsion Dixmier-Douady classes for various cases where
we relax the hypotheses on the fibre and  base.  


\section{Bundle gerbes}
We review here quickly the basic results on bundle gerbes needed to understand the proof and later examples.
The reader is referred to \cite{Mur96Bundle-gerbes, Mur2010, MurSte00Bundle-gerbes:-stable} for further details and additional references.

\subsection{Basic definitions}

Let $\pi \colon Y \to M$ be a surjective submersion and denote by  $Y^{[p]}$  the $p$-fold fibre product
$$
Y^{[p]} = \{(y_1, \dots, y_p) \mid \pi(y_1) = \dots = \pi(y_p) \} \subset Y^p .
$$
For each $i=1, \dots, p+1$ define the projection  $\pi_i \colon Y^{[p+1]} \to Y^{[p]}$ to be the map that omits the $i$-th element. 

Here and elsewhere if $Q$ and $R$ are two $U(1)$ bundles we define their {\em product} $Q \otimes R$ to be  the quotient
of the fibre product of $Q$ and $R$ by the $U(1)$ action $(q, r)z = (qz, rz^{-1})$, with the induced right action 
of $U(1)$ on equivalence classes being given by $[q, r]w= [q, rw] = [qw, r]$. In other words observe that
the fibre product is a $U(1) \times U(1)$ bundle and quotient by the subgroup $\{ (z, z^{-1}) \mid z \in U(1) \}$.

In addition if $P$ is a $U(1)$ bundle  we denote by $P^*$ the $U(1)$ bundle with the same total space as $P$ but with the action of $U(1)$ changed to its inverse, thus if 
$u\in P^*$ and $z\in U(1)$ then $z$ acts on $u$ by sending it to $uz^{-1}$.  We will  refer to $P^*$ as the {\em dual} $U(1)$ bundle to $P$.  

If $L$ and $J$ are the hermitian 
line bundles associated to $P$ and $Q$ respectively then there are canonical isomorphisms between $L\otimes J$ and the hermitian line bundle associated to 
$P\otimes Q$, as well as canonical isomorphisms between the dual line bundle $L^*$ and the hermitian line bundle associated to $P^*$.  

If $Q \to Y^{[p]}$ is a $U(1)$  bundle we define a new $U(1)$ bundle $\delta(Q) \to Y^{[p+1]}$
by
$$
\delta(Q) = \pi_1^*(Q) \otimes \pi_2^*(Q)^* \otimes \pi_3^*(Q) \otimes \cdots.
$$
It is straightforward to check that $\delta(\delta(Q)) $ is canonically trivial as a $U(1)$ bundle.

We then have:
\begin{definition}
A {\em bundle gerbe} \cite{Mur96Bundle-gerbes} over $M$ is a pair $(P, Y)$ where $Y \to M$ 
is a surjective submersion and $P \to Y^{[2]}$ is a $U(1)$ bundle   satisfying the following two conditions.
\begin{enumerate}
\item There is a {\em bundle gerbe multiplication} which is a smooth isomorphism
$$
m \colon \pi_3^*(P) \otimes \pi_1^*(P) \to \pi_2^*(P)
$$
of $U(1)$ bundles over $Y^{[3]}$.    
  \item This multiplication is associative,   that is, if we let  $P_{(y_1, y_2)}$ denote the fibre of $P$ over $(y_1, y_2)$ then    the 
following diagram commutes for all $(y_1, y_2, y_3, y_4) \in Y^{[4]}$: 
\begin{equation*}
\begin{array}{ccc}
P_{(y_1, y_2) } \otimes P_{(y_2, y_3) } \otimes P_{(y_3, y_4) }  & \rightarrow  & P_{(y_1, y_3) } \otimes P_{(y_3, y_4) }  \\
\downarrow    &     & \downarrow    \\
P_{(y_1, y_2) } \otimes P_{(y_2, y_4) }  & \rightarrow &  P_{(y_1, y_4) }\\
\end{array}
\end{equation*}

\end{enumerate}
\end{definition}

It is easy to check that for every $y \in Y$ there is a unique element  $e \in P_{(y, y)}$ such that for any $p \in Y_{(y, z)} $ we have $e p = p \in Y_{(y, z)} $ and for any $q \in Y_{(x, y)}$ we have $ q e = q \in Y_{(x, y)}$. Also for any $p \in P_{(x, y)}$ there is a unique $p^{-1} \in P_{(y, x)}$ such that $p p ^{-1} = e = p^{-1} p$.  
 
\subsection{Triviality and the Dixmier-Douady class}
\label{sec:DDclass}
Bundle gerbes are higher dimensional analogues of line bundles.  Accordingly they share many of the familiar properties of line bundles: 
just as we can pullback line bundles by smooth maps, form duals and take tensor products, we can do the same 
for bundle gerbes.      

If $(P,Y)$ is a bundle gerbe over $M$ then we can form the {\em dual} bundle gerbe $(P^*,Y)$ by setting $P^*\to Y^{[2]}$ to be the dual of the $U(1)$ bundle $P$ in the sense 
described earlier.  The process of forming duals commutes with taking pullbacks and forming tensor products and so we see the bundle gerbe multiplication on $P$ 
induces a bundle gerbe multiplication on $P^*$ in a canonical way.  

If $(P,Y)$ and $(Q,X)$ are bundle gerbes over $M$ then we can form a new bundle gerbe $(P\otimes Q,Y\times_M X)$ over $M$ called 
the {\em tensor product} of $P$ and $Q$.  Here the surjective submersion 
is  the fiber product $Y\times_M X\to M$ and $P\otimes Q$ is the $U(1)$ bundle on $(Y\times_M X)^{[2]}$ whose fibre at $((y_1,x_1),(y_2,x_2))$ is given 
by 
$$ 
P_{(y_1,y_2)}\otimes Q_{(x_1,x_2)}.   
$$
The bundle gerbe multiplication on $P\otimes Q$ is defined in the obvious way, using the bundle gerbe multiplications on $P$ and $Q$.  Note that if $Y = X$ then we can form 
the tensor product bundle gerbe in a slightly different way.  We use the original surjective submersion $Y\to M$ and we define $P\otimes Q$ to be the $U(1)$ bundle 
with fiber $P_{(y_1,y_2)}\otimes Q_{(y_1,y_2)}$ at $(y_1,y_2)\in Y^{[2]}$.  The bundle gerbe multiplication is again induced from the multiplications on $P$ and $Q$.  We 
will call the bundle gerbe $(P\otimes Q,Y)$ the {\em reduced} tensor product of $P$ and $Q$.   

A bundle gerbe $(P,Y)$ over $M$ is said to be {\em trivial} if there is a $U(1)$ bundle $Q$ on $Y$ such that $P = \delta(Q)$ and the bundle gerbe multiplication on $P$ is given by  the isomorphism  
$$ 
Q_{y_1}^*\otimes Q_{y_2} \otimes Q_{y_2}^*\otimes Q_{y_3} \cong Q_{y_1}^*\otimes Q_{y_3} 
$$ 
resulting from the  canonical pairing between $Q_{y_2}$ and $Q_{y_2}^*$.  

Just as every line bundle $L$ on $M$ has a characteristic class in $H^2(M,\ZZ)$, the Chern class $c_1(L)$ of $L$, 
every bundle gerbe $(P,Y)$ over $M$ has a characteristic class in $H^3(M,\ZZ)$.  
This characteristic class is called the {\em Dixmier-\-Douady}
class and is denoted $DD(P,Y)$.      We construct it in terms of \v{C}ech cohomology as follows.     Choose a good cover $\cU = \{U_\a\}$ of $M$  \cite{BotTu} with 
sections $s_\a \colon U_\a \to Y$ of $\pi\colon Y\to M$.      Then 
$$
(s_\a, s_\b) \colon U_\a \cap U_\b \to Y^{[2]}
$$
is a section.      Choose a section $\sigma_{\a\b} $ of $P_{\a\b} = (s_\a, s_\b)^*(P)$. That is, $\sigma_{\a\b}$ is a map such that  
$$
\sigma_{\a\b} \colon  U_\a \cap U_\b \to P
$$
with $\sigma_{\a\b}(x) \in P_{(s_\a(x), s_\b(x))}$.      Over triple overlaps we have 
$$
m(\sigma_{\a\b}(x), \sigma_{\b\c}(x)) = g_{\a\b\c}(x) \sigma_{\a\c}(x) \in P_{(s_\a(x), s_\c(x))}
$$
for  $g_{\a\b\c} \colon U_\a \cap U_\b \cap U_\c \to U(1)$.       This defines a 
co-cycle  which represents the Dixmier-Douady class
$$
DD(P, Y) = [g_{\a\b\c}] \in H^2(M, U(1)) = H^3(M, \ZZ).
$$
The Dixmier-Douady class of $P$ is the obstruction to $(P,Y)$ being trivial in the sense that $DD(P,Y)$ vanishes if and only if $(P,Y)$ is isomorphic to a trivial 
bundle gerbe.  Note also that the Dixmier-Douady class is compatible with forming tensor products in the sense that $DD(P\otimes Q,Y\times_M X) = 
DD(P,Y) + DD(Q,X)$.  Likewise for the reduced tensor product we have $DD(P\otimes Q,Y) = DD(P,Y) + DD(Q,Y)$.    

We also need to understand the image of the Dixmier-Douady class in real cohomology.  This can be defined in terms of de Rham cohomology as follows.

\subsection{Connections, curving and the real Dixmier-Douady class}

Let $\Omega^q(Y^{[p]})$ be the space of differential $q$-forms on $Y^{[p]}$.  Define a homomorphism 
\begin{equation}
\label{eq:delta-forms}
\delta \colon \Omega^q(Y^{[p]}) \to \Omega^q(Y^{[p+1]}) \quad\text{by}\quad \delta =  \sum_{i=1}^{p+1} (-1)^{i-1} \pi_i^* .
\end{equation}

These maps form the {\em fundamental  complex}
$$
0 \to \Omega^q(M) \stackrel{\pi^*}{\to} \Omega^q(Y) \stackrel{\delta}{\to}  \Omega^q(Y^{[2]}) \stackrel{\delta}{\to}
\Omega^q(Y^{[3]})    \stackrel{\delta}{\to} \dots
$$
which is exact \cite{Mur96Bundle-gerbes}.
If $(P, Y)$ is a bundle gerbe on $M$ then a {\em bundle gerbe connection} is a
connection $\nabla$ on $P$ which commutes with the bundle gerbe multiplication. If  $F_\nabla$ is the
 curvature of a bundle gerbe connection $\nabla$ then $\delta(F_\nabla) = 0$ so,  from the exactness of the fundamental complex,  $F_\nabla = \delta(f)$ for some two-form $f \in \Omega^2(Y)$. 
A choice of such an $f$ is
called a {\em curving} for $\nabla$. From the exactness of the fundamental complex we see that the 
curving is only unique up to addition of two-forms pulled back to $Y$ from $M$.  
Given a choice of curving $f$ we have  $\delta(df) = d \delta(f) =
dF_\nabla = 0$ so that $df = \pi^*(\omega)$ for a closed three-form $\omega$ on $M$ called the \emph{three-curvature}
of $\nabla$ and $f$.  The de Rham class
$$
\left[ \frac{1}{2\pi i} \omega \right] \in H^3(M, \RR)
$$
is an integral class which is the image in real cohomology of the Dixmier-Douady class of  $(P, Y)$. For convenience
let us call this the real Dixmier-Douady class of $(P, Y)$. 

\subsection{The lifting bundle gerbe}
For the sake of completeness and because we use it in the examples in the last section let us review the construction of the 
lifting bundle gerbe \cite{Mur96Bundle-gerbes}.  Let $P \to M$ be a principal $G$ bundle and note that there is a natural function $\tau \colon P^{[2]} \to G$
defined by $p_1 \tau(p_1, p_2) = p_2$. Assume moreover that $G$ has a central extension 
$$
U(1) \to \widehat G \to G.
$$
Regarding  this as a $U(1)$ bundle $\widehat G \to G$ and pulling it back with $\tau$ defines a $U(1)$-bundle $Q \to P^{[2]}$.
It is easy to check that the multiplication in $\widehat G$ induces a bundle gerbe product.  The Dixmier-Douady class of this
bundle gerbe has a well-known geometric interpretation as the obstruction to lifting the $G$ bundle $P$ to a $\widehat G$ bundle.

\section{The theorem}

\begin{theorem}
\label{thm:result} 
Let $Y \to M$ be a fibre bundle  with finite-dimensional $1$-connected fibres. Let $M  $ also be $1$-connected. 
Then any bundle gerbe $(P, Y)$ over $M$ has exact three-curvature and hence torsion Dixmier-Douady class.
\end{theorem}

As stated earlier the proof in \cite{Mur2010} is incorrect but it is  possible to fix it  as follows. Consider first the exact statement of 
the results in \cite{GotLasSni} in the case of two-forms:
 
\begin{theorem}[\cite{GotLasSni} Theorem 1]
\label{thm:one}
Let $F \to E \to B$ be a differentiable fibre bundle carrying a field $\omega$ of two-forms on the vertical bundle $\cV$, defining a closed form
on each fibre. Then there is a closed form on $E$ extending $\omega$ if and only if there is a de Rham cohomology class $c$ on $E$
whose restriction to each fibre is the class determined by $\omega$.
\end{theorem}

Although it is not spelt out in the statement the construction assumes that the class of the extension of $\omega$ is $c$.
They also prove:

\begin{theorem}[\cite{GotLasSni} Theorem 2]
\label{thm:two}
Let $F \to E \to B$ be a fibre space with $F$ and $B$ 1-connected.  If the restriction map $H^2(E, \RR) \to H^2(F, \RR)$ is not
surjective then $H^{2k}(F, \RR) \neq 0$ for all $k \geq 0$.
\end{theorem}

We use these results to prove Theorem \ref{thm:result}. Let $(P, Y)$ be the bundle gerbe. If $m \in M$ denote by $Y_m$ the fibre of $Y \to M$ over $m$. Notice that because $M$ is $1$-connected, if $m, m' \in M$  there is a unique homotopy equivalence
between $Y_m$ and $Y_{m'}$ and hence a unique identification of $H^2(Y_m)$ and $H^2(Y_{m'})$ for any choice of coefficients. 

If $m \in M$  we have a restriction map $H^2(Y, \QQ) 
\to H^2(Y_m, \QQ)$  which induces an onto map  $H^2(Y, \RR) \to H^2(Y_m, \RR)$ by Theorem \ref{thm:two}. It is easy to see that this 
implies that $H^2(Y, \QQ) \to H^2(Y_m, \QQ)$ is also onto. Indeed choose a 
basis for $H^2(Y, \QQ)$, and $H^2(Y_m, \QQ)$. Then the restriction map is given by a matrix with rational entries.  So its row reduced echelon form 
has rational entries and one can find a rational vector mapping to any rational vector.  

 Let $\nabla$ be a bundle gerbe connection with curvature
$F$ and curving $f$. Fix $y_0 \in Y_m$ and define $\iota \colon Y_m \to Y^{[2]}$ by $\iota(y) = ( y_0, y)$. Then $\pi_1 \circ \iota (y) = y$ and 
$\pi_2 \circ \iota(y) = y_0$ so that $\iota^*(F) = \iota^* \delta f = \iota^* \pi_1^* (f) - \iota^* \pi_2^*(f) = f$.  Hence $f$ restricted to 
$Y_m$ is integral and certainly rational. We deduce from Theorem \ref{thm:two} that there is a rational class in $H^2(Y, \QQ)$ extending the class defined by $f$ on any fibre 
and, moreover, it can be represented by a closed two-form $\rho$ from Theorem \ref{thm:one}.

Rationality implies that   there is some integer $n$ such that $n\rho$ is an integral two-form on $Y$. If we form
the $n$th reduced tensor  power $P^n$ of $P$ it has curving and curvature which are $n$ times the 
curving and curvature of $P$ and  we have $DD(P^n), Y) = n DD(P, Y)$.  As we are trying to show that $DD(P, Y)$
is a torsion class it suffices to show that $DD(P^n, Y)$ is a torsion class and so we may as well
assume that $n = 1$ or , in other words, that $\rho $ is integral.  

 As $Y_m$ and $M$ are 1-connected so also is $Y$ and hence $\rho$ 
defines a $U(1)$ bundle $Q \to Y$ whose curvature is $\rho$. Consider the bundle gerbe $P \otimes \delta(Q^*) \to Y^{[2]}$. This has curvature
$F - \delta(\rho)$ with curving $f - \rho$ which is zero restricted to the fibres of $Y \to M$. 
It follows that $F$ is zero restricted to 
the fibres of $Y^{[2]} \to M$ as $F = \delta(f - \rho)$.  Since the fibres of $Y\to M$ are 1-connected, the fibres of 
$Y^{[2]}\to M$ are 1-connected and so we can descend $P \otimes \delta(Q^*) $ to a 
bundle $R \to M$ by taking covariantly constant sections over the fibres of $Y^{[2]} \to M$. This descended bundle will have  connection $a$ and curvature $F_a$ whose
pull-back to $Y^{[2]}$ is the connection $\nabla$ and curvature $F - \delta(\rho)$ of $P \otimes \delta(Q^*)$. But now we have a two-form $F_a$ 
on $M$ whose pullback to $Y^{[2]}$ is zero under $\delta$. If we denote the projection from $Y^{[p]}$ to $M$ by $\pi^{[p]}$ and let $\pi_i \colon 
Y^{[p]} \to Y^{[p-1]}$ be one of the usual projections we have $\pi^{[p-1]} \circ \pi_i = \pi^{[p]}$ and in particular 
$\pi^{[2]} \circ \pi_i = \pi^{[3]}$. It follows that 
\begin{align*}
0 &= \delta ({\pi^{[2]}}^*(F_a)) \\
   & = ({\pi^{[2]}} \circ \pi_1)^*(F_a) - ({\pi^{[2]}} \circ \pi_2)^*(F_a)  + ({\pi^{[2]}}\circ \pi_3 )^*(F_a)\\
   & = ({\pi^{[3]}} )^*(F_a)
\end{align*}
and $({\pi^{[3]}})^*$ is injective so $F_a = 0$.  Hence $F - \delta(\rho) = 0$ and the bundle gerbe $P \otimes \delta(Q^*)$  has zero three-curving and thus torsion Dixmier-Douady class. But $DD(P) = DD(P \otimes \delta(Q^*))$ and thus is also torsion.  This proves Theorem \ref{thm:result}. 

\begin{note}
For the interested reader we note that the mistake in the original proof in \cite{Mur96Bundle-gerbes} was to claim that because the forms $f - \rho$ and $d(f -\rho)$
were  vertical in the sense of restricting to zero on fibres the form $f - \rho$ descended to $M$.  This is, of course, not true.  What is 
true is that if a form $\mu$ on the total space of a fibre bundle and its exterior derivative $d\mu$ are vertical in the stronger sense of vanishing when contracted with any  vertical vector, then $\mu$ descends to the base. 
\end{note}

Call a bundle gerbe $(P,Y)$ over $M$ a {\em finite bundle gerbe} if $Y$ is a  fibre bundle over $M$ with finite-dimensional fibres.  We can
restate Theorem \ref{thm:result} as:

\begin{theorem}
\label{thm:finite}
Let $(P, Y)$ be a finite bundle gerbe over $M$. If $M$ and the fibres of $Y \to M$ are $1$-connected
then $(P, Y)$ has torsion Dixmier-Douady class. 
\end{theorem}

We show now how to extend this result to the case of fibres which are only connected.  First we have

\begin{proposition}
\label{prop:three-sphere}
Let $(P, Y)$ be a finite bundle gerbe over $S^3$ with connected fibre $F$. Then $(P, Y)$ has torsion and hence zero Dixmier-Douady class. 
\end{proposition}

\begin{proof}

Form the universal cover $\widetilde{Y}\to Y$.  Then we have a 
diagram 
$$ 
\xy
(-10,0)*+{\widetilde{Y}}="1";
(10,0)*+{Y}="2";
(0,-12)*+{S^3}="3";
{\ar_p "1";"2"};
{\ar_-{\tilde{\pi}} "1";"3"}; 
{\ar^-{\pi} "2";"3"};
\endxy 
$$
Since $Y\to S^3$ is locally trivial and $\widetilde{Y}\to Y$ is a covering space, $\widetilde{Y}\to S^3$ is locally trivial with fiber $\widetilde{F}$, where $\widetilde{F}$ denotes the pullback of 
$\widetilde{Y}\to Y$ under the inclusion of the fiber $F\subset Y$.  Consider the long exact homotopy sequences of the fibrations $\widetilde{Y}\to S^3$ and 
$Y\to S^3$, by naturality we have a commutative diagram 
$$ 
\xymatrix{ 
\cdots \ar[r] & \pi_2(S^3) \ar[d] \ar[r] & \pi_1(\widetilde{F}) \ar[d] \ar[r] & \pi_1(\widetilde{Y}) \ar[d] \ar[r] & \pi_1(S^3) \ar[d] \\ 
\cdots \ar[r] & \pi_2(S^3) \ar[r] & \pi_1(F) \ar[r] & \pi_1(Y) \ar[r] & \pi_1(S^3) } 
$$ 
from which we conclude that $\pi_1(\widetilde{F}) = 1$.  Note that since $\widetilde{Y}\to Y$ is a covering space and covering spaces pullback to covering spaces, 
$\widetilde{F}\to F$ is also a covering space, in fact $\widetilde{F}$ is the 
universal covering of $F$. 

In such a case as this  the map $p \colon \widetilde{Y} \to    Y$ induces a map $p^{[2]} \colon \widetilde{Y}^{[2]} \to Y^{[2]}$
and we can pull back $P \to Y^{[2]}$ to form a bundle gerbe $((p^{[2]})^*(P), \widetilde{Y})$. It is straightforward
from the explicit construction of the Dixmier-Douady class in Section \ref{sec:DDclass} to show that $DD(((p^{[2]})^*(P), \widetilde{Y}) = DD(P, Y)$. But as $\tilde F$ is $1$-connected Theorem \ref{thm:result} gives us that $DD(P, Y)$
is torsion.  But $H^3(S^3, \ZZ) = \ZZ$ so that  $DD(P, Y) = 0$.
\end{proof}

We now have 

\begin{theorem}
\label{thm:improved}
Let $(P, Y)$ be a finite bundle gerbe over a simply connected manifold $M$ with connected fibre $F$. Then $(P, Y)$ has torsion Dixmier-Douady class. 
\end{theorem}
\begin{proof}
As $M$ is a simply connected manifold Hurewicz's theorem implies that the Hurewicz homomorphism
$$
h \colon \pi_3(M) \to H_3(M, \ZZ)
$$
 is onto. Recall that $h $ is defined by choosing a generator $e \in H_3(S^3, \ZZ)$ and letting $h([\alpha]) = 
 \alpha_*(e)$.  Recall also that there is a homomorphism 
$$
 \rho \colon H^3(M, \ZZ) \to \Hom(H_3(M, \ZZ) , \ZZ)
$$
defined by pairing the cohomology and homology classes whose kernel is the torsion subgroup of $H^3(M, \ZZ)$. 

Let $(P, Y)$ be a finite bundle gerbe with Dixmier-Douady class $DD(P, Y)$.  Then $\rho(DD(P, Y))$
can be determined by evaluating it on classes of the form $h([\alpha])$ to get 
$\rho(DD(P, Y)) = \rho(DD(P, Y))(h([\alpha])) = \a^*(DD(P, Y)) = DD(\a^*(P), \a^*(Y)) = 0$.
Hence $DD(P, Y)$ is torsion.   
\end{proof}

\section{Examples}

We consider some examples to see what can be said about the necessity of the conditions in Theorem \ref{thm:improved}.
Note first that a bundle gerbe over $M$ restricts to a bundle gerbe over any connected component of $M$ so there is nothing
of interest to be lost by assuming that $M$ is connected which we do henceforth. 
Before considering the constructions we need to make two general remarks. First we introduce some notation: 
given a map $a \colon Y^{[p]} \to A$ for some abelian group $A$ we define $ \delta(a) \colon Y^{[p+1]} \to A$ by 
$$
\delta(a) = \pi_1^*a - \pi_2^*a + \pi_3^*a - + \cdots  
$$
Secondly if $Y \to M$  is a surjective submersion, then one way to define a bundle gerbe is to consider a function 
$c \colon Y^{[3]} \to  U(1)$, take $P \to Y^{[2]}$ to be the 
trivial $U(1)$ bundle, and define a bundle gerbe product by 
$$
((y_1, y_2), \a)((y_2, y_3), \b) = ((y_1, y_3), c(y_1, y_2, y_3)  \a\b).
$$
This product is associative if and only if $\delta(c) = 1$.
A bundle gerbe connection for $Q$ is a  one-form $A$ on $Y^{[2]}$ satisfying $\delta(A) = h^{-1} dh$ where $\delta(A)$ 
is defined in equation \eqref{eq:delta-forms}. 
The curvature of $A$ is $dA$. We can then define curving and  three-curvature in the usual way.

\subsection{Bundle gerbes from open covers}
It is perhaps worth remarking that if $[g_{\a\b\c}] $ is a representative co-cycle for a class in $H^3(M, \ZZ)$
with respect to an open cover $\cU = \{ U_\a \}_{\a \in I}$ of $M$ we can define $Y_\cU$ to be the 
disjoint union of the open sets in the cover $\cU$ with the obvious projection $Y_{\cU} \to M$. Then $g_{\a\b\c}$
defines a function $c \colon Y^{[3]}_{\cU} \to U(1)$ as above and it is easy to check that this defines a finite
bundle gerbe with Dixmier-Douady class $[g_{\a\b\c}]$.  In this example $Y_{\cU} \to M$ is a surjective submersion but is, of course, unlikely to be a fibre bundle. 

\subsection{Cup product bundle gerbes.} A nice way to construct examples of bundle gerbes is via the cup product construction (see for example \cite{Bry,Stu}).  
Suppose we are given geometric representatives of classes $\a$ in 
$H^2(M,\ZZ)$ and $\b$ in $H^1(M,\ZZ)$ corresponding to a principal $U(1)$ bundle $Q$ on $M$ and a smooth map $f\colon M\to S^1$ respectively.  Then there is a 
bundle gerbe over $M$ with Dixmier-Douady class equal to the cup product $\a\cup \b$. There are two ways in which this can 
be described which are of interest to us. In the first case the bundle gerbe is of the form $(P, Y)$ where $Y$ is the 
$\ZZ$-bundle $f^*\RR$ and where $\RR \to S^1$ is the universal $\ZZ$ bundle.  In the second case it is of 
the form $(P, Y)$ where $Y$ is the $\ZZ \times U(1)$ fibre bundle
which is the fibre product $f^*\RR \times_M Q$ of $f^*(\RR) $ and $Q$. Notice that in both cases $Y$ is disconnected.

Let us consider the first case In more detail.  Take $Y  = f^*\RR$.  Then there is a map $\tau\colon Y^{[2]}\to \ZZ$ defined by 
$y_2 = y_1\tau(y_1,y_2)$ for $(y_1,y_2)\in Y^{[2]}$ and we can define $P\to Y^{[2]}$ to be the $U(1)$ bundle with fibre at $(y_1,y_2)$ given by $Q_m^{\otimes \tau(y_1,y_2)}$, 
where $\pi(y_1) = \pi(y_2) = m$.  

Likewise, in the second case $Y$ is the fibre product $Q\times_M f^*\RR$ so that $Y$ is a principal $U(1)\times \ZZ$ bundle over $M$.  The group 
$U(1)\times \ZZ$ fits into a central extension of Lie groups 
$$ 
U(1)\to U(1)\times \ZZ\times U(1)\to U(1)\times \ZZ 
$$ 
where the product on $U(1)\times \ZZ\times U(1)$ is defined by 
$$ 
(z_1,n_1,w_1)\cdot (z_2,n_2,w_2) = (z_1z_2, n_1 + n_2, w_1w_2z_1^{n_1}) 
$$ 
We refer the reader to \cite{Bry} for more details.  The bundle gerbe $(P,Y)$ is then given by the lifting bundle gerbe construction.  One can check (see Corollary 4.1.15 of \cite{Bry}) that the 
Dixmier-Douady class $DD(P,Y)$ is given by the cup product $\a\cup \b$.

As an example let us take $M = S^2\times S^1$.  We let $\a$ denote the class in $H^2(M,\ZZ)$ defined by pulling back the Hopf bundle $S^3\to S^2$ 
via the projection to $S^2$ and we let $\b$ denote the class in $H^1(M,\ZZ)$ defined by the projection to $S^1$.  In the first construction we take $Y = S^2 \times \RR \to S^2 \times S^1$ and in the second  $Y = S^3\times \RR \to S^2 \times S^1$.

This example can be greatly generalized.  Suppose that $G$ is a compact, simple 1-connected Lie group with maximal torus $T$.  Let $\mathfrak{t}$ denotes the 
Lie algebra of $T$.  Then there is a natural principal bundle $G\times \mathfrak{t}\to G/T\times T$ with structure group $T\times \pi_1(T)$.  Given a bilinear form $b$ on the Lie algebra $\mathfrak{t}$ of 
$T$, one can define a central extension of groups $U(1)\to T\times \pi_1(T)\times U(1)\to T\times \pi_1(T)$ (see for example \cite{PS}) and so we can form the corresponding 
lifting bundle gerbe.       

\subsection{Bundle gerbes on unitary groups.} 
 Theorem~\ref{thm:improved} implies that  a finite bundle gerbe with connected fibres over a simply connected, simple compact Lie group $G$
must be torsion.  In particular the basic bundle gerbe corresponding to the standard  generator of $H^3(G, \ZZ)$ cannot be 
a finite bundle gerbe with connected fibres.  We have shown in \cite{MurSte} that when $G = SU(n)$ it is possible to realise the bundle gerbe with 
Dixmier-Douady class the standard generator of $H^3(SU(n), \ZZ)$ as a finite bundle gerbe with disconnected fibres as follows
We define
$$
Y = \{ (X, \lambda) \mid \det(X - \lambda I) \neq 0 \} \subset SU(n) \times Z,
$$
where $Z$ denotes the {\em set} $U(1)$ with the identity element removed.  
A point in $Y^{[2]}$ can be thought of as a triple $(X, \a, \b)$ where neither of $\a$ or $\b$ is an eigenvalue of $X$. We define a 
hermitian line bundle over $Y^{[2]}$ by taking the fibre at $(X,\a, \b)$ to be the determinant of the sum of the eigenspaces 
of $X$ lying between $\a$ and $\b$ on $Z$, with respect to a certain ordering on $Z$.  The corresponding $U(1)$ bundle is the required bundle gerbe.  Of course in this case
the fibres of $Y \to SU(n)$ are disconnected and it is not, in fact, a fibre bundle.  

Other constructions of the basic bundle gerbe on a compact Lie group with the fibres of $Y$ either disconnected or infinite-dimensional have been considered by other
authors and are reviewed in the introduction to \cite{MurSte}.

\subsection{A bundle gerbe on the three-torus} 
Consider $T^3 = S^1 \times S^1 \times S^1$ and 
the projection $Y = \RR^3 \stackrel{\pi}{\to} T^3$ induced by the standard projection $\RR \to S^1$ given by  $t \mapsto \exp(2\pi i t)$. 
Notice that the fibres of $\pi\colon Y \to T^3$ are disconnected and the base is, of course, not simply connected. Using constructions from 
\cite{Stu} we show how to 
construct the bundle gerbe whose Dixmier-Douady class is the natural generator of $H^3(T^3, \ZZ)$.  

We will write $x = (x^1,x^2,x^3)$ for a vector $x$ in $Y$.  
Note that  $(x, y) \in Y^{[2]}$ if and only if $x - y \in \ZZ^3$. 
If $(x, y, z) \in Y^{[3]} \subset  \RR^3 \times \RR^3 \times \RR^3 $ define
  $\gamma(x, y, z) =   (y^1-z^1)(x^2-y^2)x^3$ and  $c(x, y, z) = \exp(2 \pi i \gamma(x, y, z))$.
If $(x, y, z, w) \in Y^{[4]}$ then 
$$
\delta(\gamma)(x, y, z, w) = \gamma(y-x, z-y, w-z)  \in \ZZ
$$
so that $\delta(c) = 1$ and this defines a bundle gerbe. Writing points of $Y^{[3]}$ as $(x, y, z)$ we can 
denote the projections as $x$, $y$, $z$ and we have $\RR^3$ valued differential forms $dx$, $dy$ and $dz$. As $x$, $y$ and $z$ differ by integers these forms are all equal and we will denote the resulting form by $\theta = (\theta^1, \theta^2, \theta^3)$. We can 
similarly define an $\RR^3$ valued one-form on $Y$ and pulling it back by either of the projections $Y^{[2]} \to Y$ gives the form $\theta$ so we will denote the form on $Y$ as $\theta$ as well. Finally notice that $\theta^i$, on $Y$, is the pull-back from the $i$th copy 
of $U(1)$ of the one-form $d\theta/{2\pi}$ which has total integral one on $U(1)$. It is now easy to check that $A = -2\pi i (x^1 - y^1)(x^2)\theta^3$ and $f = 2\pi i x^1 
\theta^2 \wedge \theta^3$ give a connection and curving for this bundle gerbe. Notice that the curvature is 
$2\pi i\, \theta^1 \wedge \theta^2 \wedge \theta^3$ so the real Dixmier-Douady class is $(1/8\pi^3) d\theta^1 \wedge d\theta^2 \wedge d\theta^3$ on $T^3$ as we require.


\end{document}